\documentclass[12pt]{article}
\usepackage{graphics}
\usepackage{graphicx}
\usepackage{color}

\usepackage{amsfonts}
\usepackage{eufrak}

\usepackage{amsmath}
\usepackage{amstext}
\usepackage{amsopn}
\usepackage{amsbsy}
\usepackage{amscd}
\usepackage{amsxtra}
\usepackage{amsthm}

\numberwithin{equation}{section}

\newcommand{\beq}{\begin{equation}}
\newcommand{\eeq}{\end{equation}}
\newcommand{\pp}{\partial}
\newcommand{\R}{{\mathbb R}}

\newcommand{\ds}{\displaystyle}

 \title{{\scshape\normalsize Mathematics Division, National Center for
         Theoretical Sciences at Taipei}\\
         {\scshape\large NCTS/TPE-Math Technical Report 2007-004} \\~\\
         Half-Skyrmions and Spike-Vortex Solutions of Two-Component Nonlinear Schr\"{o}dinger Systems}

\author{Tai-Chia Lin  \thanks{ Department of Mathematics, National Taiwan
University, Taipei, Taiwan 106. email : tclin@math.ntu.edu.tw}
\thanks{National Center for Theoretical Sciences,Taipei Office }\ and Juncheng Wei
\thanks{Department of Mathematics, The Chinese University of Hong Kong,
Shatin, Hong Kong. email: wei@math.cuhk.edu.hk} }

\date{}

\newtheorem{thm}{\large T\small HEOREM\large}[section]
\newtheorem{lem}[thm]{\large L\small EMMA\large}

\theoremstyle{remark}

\begin{document}
\maketitle

\begin{abstract}

\ \ \ Recently, skyrmions with integer topological charges have
been observed numerically but have not yet been shown rigorously
on two-component systems of nonlinear Schr\"{o}dinger equations
(NLSE) describing a binary mixture of Bose-Einstein condensates
(cf.~\cite{BCS} and \cite{SR}). Besides, {\it half-skyrmions}
characterized by half-integer topological charges can also be
found in the nonlinear $\sigma$ model which is a model of the
Bose-Einstein condensate of the Schwinger bosons (cf.~\cite{M}).
Here we prove rigorously the existence of {\it half-skyrmions}
which may come from a new type of soliton solutions called
spike-vortex solutions of two-component systems of NLSE on the
entire plane $\mathbb{R}^2$. These spike-vortex solutions having
spikes in one component and a vortex in the other component may
form half-skyrmions. By Liapunov-Schmidt reduction process, we may
find spike-vortex solutions of two-component systems of NLSE.

\end{abstract}

\section{Introduction}
\ \ \ \ Spikes and vortices are important phenomena in
one-component nonlinear Schr\"{o}dinger equations (NLSE) having
applications in many physical problems, especially in
Bose-Einstein condensation. In the last two decades, there have
been many analytical works on both spikes and vortices,
respectively. One may refer to~\cite{n} for a good survey on
spikes, and \cite{BBH}, \cite{L1} and \cite{PR} for survey on
vortices. Recently, a double condensate i.e. a binary mixture of
Bose-Einstein condensates in two different hyperfine states has
been observed and described by two-component systems of NLSE
(cf.~\cite{PS}). It would be possible to find spike and vortex
solutions from two-component systems of NLSE. However, until now,
there is no result to deal with spike-vortex solutions having
spikes in one component and a vortex in the other component. In
this paper, we want to find such solutions and investigate the
interaction of spikes and vortices.

Ordinary skyrmions being of topological solitons resemble
polyhedral shells which look like closed loops, possibly linked or
knotted. Besides their intrinsic fundamental interest, skyrmions
have important applications in nuclear physics~(cf.~\cite{1}), and
analogous structures are postulated for early universe
cosmology~(cf.~\cite{C}). To get a skyrmion, one may consider the
multi-component wave function which may introduce the extra
internal degrees of freedom and result in a nontrivial structure
characterized by topological charges. For two-dimensional Skyrme
model, skyrmions have been investigated by the method of
concentration-compactness~(cf.~~\cite{LY}). In a double
condensate, skyrmions with integer topological charges have been
observed by numerical simulations on two-component systems of NLSE
(cf.~\cite{BCS} and \cite{SR}). Besides, {\it half-skyrmions}
characterized by half-integer topological charges can also be
found in the nonlinear $\sigma$ model which may describe the
Bose-Einstein condensate of the Schwinger bosons (cf.~\cite{M}).
Here we prove rigorously the existence of {\it half-skyrmions} in
a double condensate using spike-vortex solutions of two-component
systems of NLSE.

To get spike-vortex solutions, we study soliton solutions of
time-dependent two-component systems of NLSE as follows:
\begin{equation}\label{id1_1}
-\sqrt{-1}\,\frac{\partial \psi_j}{\partial
t}=\triangle\psi_j+\sum_{i=1}^2\,\beta_{ij}|\psi_i|^2\psi_j\,,\quad\hbox{
for }\: x\in\mathbb{R}^n\,, t>0\,, j=1, 2\,,
\end{equation} where the spatial dimension $n=2$ and
$\psi_j=\psi_j(x,t)\in\mathbb{C}$ for $j=1, 2$. The system
(\ref{id1_1}) is a standard model to describe a double condensate.
Physically, $\psi_j$'s are the corresponding complex-valued wave
functions, and the coefficients $\beta_{ij}\sim -a_{ij}$ for
$i,j=1, 2$, where $a_{jj}$'s and $a_{12}=a_{21}$ are the
intraspecies and interspecies scattering lengths, respectively.
When the spatial dimension is one, i.e. $n=1$, it is well-known
that the system (\ref{id1_1}) is integrable, and there are many
analytical and numerical results on soliton solutions of coupled
nonlinear Schr\"{o}dinger equations (e.g.~\cite{H}, \cite{HS},
\cite{KL}). Recently, from physical experiments (cf.~\cite{CTW}),
even three-dimensional solitons have been observed in
Bose-Einstein condensates. It is natural to believe that there are
two-dimensional (i.e.~$n=2$) solitons in double condensates.
However, when the spatial dimension is two, the system
(\ref{id1_1}) becomes non-integrable and has only few results on
two-dimensional solitons. This may lead us to study
two-dimensional soliton solutions of the system (\ref{id1_1}) and
find different types of solitons.

In the vicinity of a Feshbach resonance, scattering lengths
$a_{ij}$'s depend sensitively on the magnitude of an externally
applied magnetic field (cf.~\cite{3} and \cite{4}), allowing the
magnitude and sign of $\beta_{ij}$'s to be tuned to any value.
Generically, when both $\beta_{jj}$'s are positive, the
system~(\ref{id1_1}) is of self-focusing and has bright soliton
solutions on the associated two components. On the other hand,
when both $\beta_{jj}$'s are negative, the system~(\ref{id1_1}) is
of self-defocusing and has dark soliton solutions on the
associated two components. Here we have interest on the case that
$\beta_{11}$ and $\beta_{22}$ have opposite signs which may result
in a new type of soliton solutions called spike-vortex solutions
of the system~(\ref{id1_1}) i.e. one component has spikes and the
other component has a vortex. Furthermore, we may obtain
half-skyrmions by these spike-vortex solutions.

To obtain soliton solutions of the system (\ref{id1_1}), we set
\begin{equation}\label{id1_2}
\psi_j(x,t)=e^{\sqrt{-1}\,\lambda_jt}\cdot u_j(x)\,,\quad
u_j\in\mathbb{C}\,,\: j=1, 2\,.
\end{equation}
Substituting (\ref{id1_2}) into (\ref{id1_1}), we may obtain a
two-component system of semilinear elliptic partial differential
equations given by
\begin{equation}
\begin{cases}\label{id1_3}
\triangle u_1-\lambda_1u_1+\beta_{11}\,|u_1|^2\,u_1+\beta_{12}\,u_1\,|u_2|^2=0,\\
\triangle
u_2-\lambda_2u_2+\beta_{22}\,|u_2|^2\,u_2+\beta_{12}|u_1|^2u_2=0,
\end{cases}
\end{equation}
where $\beta_{12}$ is the coupling constant.
In~\cite{LW1}-\cite{LW2}, we studied the ground state solutions of
(\ref{id1_3}) for the case that $\beta_{11}>0$, $\beta_{22}>0$ and
$u_j$'s are positive functions. Namely, we investigated the
following problem:
\begin{equation}\label{id1_4}
\begin{cases}
\triangle u-\lambda_1u+\beta_{11}u^3+\beta_{12}uv^2=0,\\
\triangle v-\lambda_2v+\beta_{22}v^3+\beta_{12}u^2v=0,\\
u,v>0,\ u,v\in H^1(\mathbb{R}^n)\,.
\end{cases}
\end{equation}
Due to each $\lambda_j>0$ and $\beta_{jj}>0$, both $u$ and $v$
components have attractive self-interaction which may let spikes
occur in these two components. One the other hand, when $n=2$,
$\lambda_j<0$ and $\beta_{jj}<0, j=1, 2$, and $u_j$'s are
complex-valued solutions of (\ref{id1_3}), vortices may exist in
both $u_1$ and $u_2$ components (cf.~\cite{ll}).

In this paper, we study the case that $n=2$, $\lambda_1,
\beta_{11}>0$, $\lambda_2, \beta_{22}<0$, and $u_1$ is positive
but $u_2$ is complex-valued function. Without loss of generality,
we may assume that
$$\lambda_1=-\lambda_2=\beta_{11}=-\beta_{22}=1\,.$$ Namely, we
study the following system:
\begin{equation}\label{id1_5}
\begin{cases}
\triangle u-u+u^3+\beta u|v|^2=0~~\quad\hbox{ in }\:\mathbb{R}^2\,,\\
\triangle v+v-|v|^2v+\beta u^2v=0\quad\hbox{ in }\:\mathbb{R}^2\,,
\end{cases}
\end{equation}
where $u>0$ and $v\in\mathbb{C}$. To get skyrmions, a defining
property of the skyrmion is that the atomic field approaches a
constant value at spatial infinity (cf.~\cite{WAS}). Hence we may
set boundary conditions of the system~(\ref{id1_5}) as follows:
$u(x)\to 0$ and $|v(x)|\to 1$ as $|x|\to\infty$.

As~$\beta=0$, the first equation of the system~(\ref{id1_5})
becomes a standard nonlinear Schr\"{o}dinger equation given by
\begin{equation}\label{id1_6}
\triangle u-u+u^3=0,\quad u\in H^1(\mathbb{R}^2)
\end{equation}
which has a unique least-energy spike solution $w=w(r), r=|x|$
satisfying $w'(r)<0$ for $r>0$ and
\begin{equation}\label{id1_7}
w(r)=A_0r^{-\frac{1}{2}}e^{-r}\Big[1+O\Big(\frac{1}{r}\Big)\Big],
\quad
w'(r)=-A_0r^{-\frac{1}{2}}e^{-r}\Big[1+O\Big(\frac{1}{r}\Big)\Big].
\end{equation} Actually, this is also a typical spike solution of
nonlinear Schr\"{o}dinger equations. On the other hand, as
$\beta=0$, the second equation of the system~(\ref{id1_5}) can be
written as
\begin{equation}\label{id1_8}
\triangle v+v-|v|^2v=0,\quad v=v(z)\in \mathbb{C}\quad\hbox{ for
}\:z\in\mathbb{C}\cong\mathbb{R}^2\,,
\end{equation}
which is of conventional Ginzburg-Landau equations
(cf.~\cite{BBH}) having a symmetric vortex solution of degree
$d\in \mathbb{Z}\backslash\{0\}$ with the following form
\begin{equation}
\label{vd}
v_d(z)=S_d(r)e^{\sqrt{-1}\,d\theta},
\end{equation}
 where $S_d(r)$ satisfies
\begin{equation}\label{id1_15}
\begin{cases}
S_d''+\frac{1}{r}S_d'-\frac{d^2}{r^2}S_d+S_d-S_d^3=0,\\
S_d(0)=0,\ S_d(+\infty)=1,
\end{cases}
\end{equation}
and
\begin{equation}\label{id1_9}
S_d'(r)>0,\quad S_d
(r)=1-\frac{d^2}{2r^2}+O\Big(\frac{1}{r^4}\Big)\quad \mathrm{as}\
r\rightarrow+\infty\,.
\end{equation} Here we want to prove that when $\beta$ increases or decreases,
 there exist spike-vortex solutions of the
system~(\ref{id1_5}). This may become the first paper to
illustrate such solutions of two-component nonlinear
Schr\"{o}dinger systems.

The main purpose of this paper is to construct a spike-vortex
solution $(u,v)$ of the system~(\ref{id1_5}) such that $u\sim w$
and $v\sim v_d$. Actually, the main difficulty of this paper is to
study the interaction of typical spike and vortex solutions. Our
first result shows that even a weak repulsive interaction
($\beta>0$ being small) can produce abundant bound states by
solving the system~(\ref{id1_5}). More precisely, we have
\begin{thm}\label{thm1_1}
Let $n=2$, $d\in\mathbb{N}$ and $k\geq2$ satisfy \begin{itemize}
\item[(i)]~~$k\geq 2$ is any positive integer if $d=1$,
\item[(ii)]~$2(d-1)\not\equiv 0 \mod k$ i.e.~there does not exist
any integer $\mu$ such that $2(d-1)=k\mu$ if $d\geq
2$.\end{itemize} Then for $\beta>0$ sufficiently small, the
problem (\ref{id1_5}) has a solution $(u_\beta,v_\beta)$
satisfying $u_\beta(z)>0$ for $z\in\mathbb{C}$, $v_\beta(0)=0$
with degree $d$, and $u_\beta(z)\to 0\,, |v_\beta(z)|\to 1$ as
$|z|\to\infty$. Moreover, as $\beta\to 0+$, $(u_\beta,v_\beta)$
has the following asymptotic form
\begin{equation}\label{id1_12}
\begin{cases}
u_\beta(z)=\ds\sum_{j=1}^kw\big(z-\xi_j^\beta\big)+O(|\beta|),\\
v_\beta(z)=S_d(r)e^{\sqrt{-1}\,d\theta} e^{ \sqrt{-1} \psi_\beta
(z)}, \ \psi_\beta (z) =O(|\beta|)\in \mathbb{C},
\end{cases}
\end{equation}
where $w$ is the unique radial solution of (\ref{id1_6}),  $\big<\xi_1^\beta,...,\xi_k^\beta\big>$ forms a regular k-polygon with
\begin{equation}\label{id1_13}
\xi_j^\beta=l_\beta e^{\sqrt{-1}\frac{2\pi(j-1)}{k}},\ j=1,...,k,
\end{equation}
and $l_\beta\to +\infty$ as $\beta\to 0+$ satisfying
$l_\beta=\hat{l}_\beta+O(1)$, where $\hat{l}_\beta$ satisfies
\begin{equation}\label{id1_14}
\hat{l}_\beta^{\frac{5}{2}}e^{-2\hat{l}_\beta\sin\frac{\pi}{k}}=\beta.
\end{equation}
\end{thm}
A picture of $(u_\beta,v_\beta)$ is given by\\

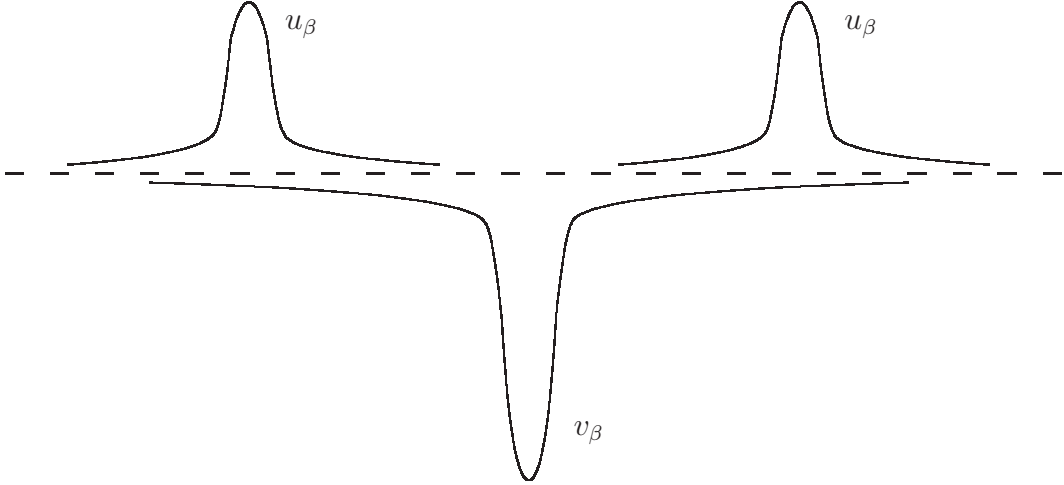
\begin{center}
\setlength{\unitlength}{1.2cm}
\begin{picture}(6,6)(-3,-4)
%%%%%%%%%%%%%%%%%%%%%%%%%%%%%%%%%%%%%%%%%%%%%%%%%%%
\multiput(-6.8,-.1)(1,0){12}
{\line(1,0){0.2}}
\multiput(-6.3,-.1)(1,0){12}
{\line(1,0){0.2}}
%%%%%%%%%%%%%%%%%%%%%%%%%%%%%%%%%%%%%%%%%%%%%%%%%%%
\put(-3.7,1.5){$u_\beta$}
\put(2.5,1.5){$u_\beta$}
\put(-.5,-3){$v_\beta$}
%%%%%%%%%%%%%%%%%%%%%%%%%%%%%%%%%%%%%%%%%%%%%%%%%%%
\qbezier(0,0)(1.45,0.1)(1.62,0.35)
\qbezier(1.62,0.35)(1.7,0.45)(1.8,1.4)
\qbezier(1.8,1.4)(1.9,1.8)(2,1.8)
\qbezier(2,1.8)(2.1,1.8)(2.2,1.4)
\qbezier(2.2,1.4)(2.3,0.45)(2.375,0.35)
\qbezier(2.375,0.35)(2.45,0.1)(4.1,0)
%%%%%%%%%%%%%%%%%%%%%%%%%%%%%%%%%%%%%%%%%%%%%%%%%%%
\qbezier(-6.1,0)(-4.65,0.1)(-4.48,0.35)
\qbezier(-4.48,0.35)(-4.4,0.45)(-4.3,1.4)
\qbezier(-4.3,1.4)(-4.2,1.8)(-4.1,1.8)
\qbezier(-4.1,1.8)(-4,1.8)(-3.9,1.4)
\qbezier(-3.9,1.4)(-3.8,0.45)(-3.715,0.35)
\qbezier(-3.715,0.35)(-3.65,0.1)(-2,0)
%%%%%%%%%%%%%%%%%%%%%%%%%%%%%%%%%%%%%%%%%%%%%%%%%%%
\qbezier(-5.2,-.2)(-1.65,-0.3)(-1.48,-0.65)
\qbezier(-1.48,-0.65)(-1.4,-0.75)(-1.3,-1.7)
\qbezier(-1.3,-1.7)(-1.2,-3.49)(-1,-3.5)
\qbezier(-1,-3.5)(-0.8,-3.49)(-0.7,-1.7)
\qbezier(-0.7,-1.7)(-0.6,-0.75)(-0.525,-0.65)
\qbezier(-0.525,-0.65)(-0.45,-0.3)(3.2,-.2)
\end{picture}
\end{center}

In \cite{LW1}-\cite{LW2}, the positive sign of $\beta$ may give
inter-component {\bf attraction} when the inter-component
interaction is only for spikes. Conversely, the positive sign of
$\beta$ may contribute inter-component {\bf repulsion} when the
inter-component interaction is for spikes and a vortex. The
inter-component repulsion between $u$ and $v$ components may
balance with self-attraction in $u$-component so a new kind of
soliton solutions called spike-vortex solutions of two-component
nonlinear Schr\"{o}dinger systems can be obtained in
Theorem~\ref{thm1_1}.

For $\beta<0$, we may consider the radial solution of
(\ref{id1_5}) given by \beq\label{ra1} u=u(r)\,,\quad
v=f(r)\,e^{\sqrt{-1}d\theta}\,, \eeq
 satisfying \beq\label{ra2}
\left\{\begin{array}{lll}
u^{\prime\prime}+\frac{1}{r}\,u^\prime-u+u^3\,+\beta f^2\,u
&=0\,,&\forall r>0\,, \\
f^{\prime\prime}+\frac{1}{r}\,f^\prime-\frac{d^2}{r^2}\,f+f-f^3+\beta
u^2\,f &=0\,, &\forall r>0\,,\\
u^\prime(0)=0, u(+\infty)=0, & & \\ f(0)=0, f(+\infty)=1\,, & &
\end{array}\right.\eeq where $d\in\mathbb{N}$ and $(r,\theta)$ is
the polar coordinates in $\mathbb{R}^2$. Then we have the
following existence theorem:

\begin{thm}
\label{thm1_2} Assume that $n=2$ and $\beta<0$. Then the problem
(\ref{id1_5}) has radially symmetric solutions of the following
form: \beq\label{ra11} u=u(r)\,,\quad v=f(r)\,e^{\sqrt{-1}d\theta}
\eeq where $ u(r)$ is strictly decreasing,  $ f(r)$ is strictly
increasing and $(u, f)$ satisfies (\ref{ra2}).
\end{thm}

\noindent {\bf Remarks:} \\
\noindent 1.~It is easy to see that the solution $(u(r), f(r))$ is
unique for the system~(\ref{ra2}) if $\beta$ is small enough. It
is an interesting question to study the uniqueness for general
$\beta<0$.

\noindent 2.~Note that the solution $(u(r), f(r) e^{\sqrt{-1} d
\theta})$ is not a global minimizer for the corresponding energy
functional of (\ref{id1_5}) since the equation of $u$ is
superlinear. It is conjectured that there exists a least energy
solution (see definitions in \cite{LW2}) satisfying $E[u(r), f(r)
e^{\sqrt{-1} d \theta}]\leq E[u, v]$ for any  solution $(u, v)$ of
(\ref{id1_5}) with $ \mbox{deg} (v)=d$. Here the corresponding
energy functional is defined as
\begin{equation}
E[u, v]= \frac{1}{2} \int_{\R^2} (|\nabla u|^2 +u^2 )-\frac{1}{4} \int_{\R^2} u^4 +\frac{1}{2} \int_{\R^2} |\nabla v|^2 +\frac{1}{4} \int_{\R^2} (1-|v|^2)^2 -\frac{\beta}{2} \int_{\R^2}  u^2 |v|^2.
\end{equation}

To get half-skyrmions, we may define a $S^2$-valued map
\beq\label{ndef} \overrightarrow{n}\equiv\frac{1}{\sqrt{u^2+v_1^2+v_2^2}}\,\left(\begin{array}{ll}&v_1\\
&v_2\\ &u\end{array}\right)\,, \eeq where $(u,v_1+\sqrt{-1}\,v_2)$
is the spike-vortex solution obtained in Theorem~\ref{thm1_1} and
\ref{thm1_2}. Generically, the topological charge of $S^2$-valued
maps is defined by (cf.~\cite{DFN}) \beq\label{topo}
Q=\frac{1}{4\pi}\,\int_{\mathbb{R}^2}\overrightarrow{n}\cdot\left(\pp_x
\overrightarrow{n}\wedge\pp_y\overrightarrow{n}\right)\,dx\,dy\,.
\eeq Then we have
\begin{thm}
\label{thm1_3} The $S^2$-map defined by (\ref{ndef}) is of
half-skyrmions with topological charge~$\frac{d}{2}$.
\end{thm}

Throughout the rest of the paper, we assume that
\begin{equation}\label{ldef}
\hat{l}_\beta-\gamma<l<\hat{l}_\beta+\gamma
\end{equation}
for some suitable $\gamma$. Note that
\begin{equation}\label{id1_17}
\hat{l}_\beta=\frac{1}{2\sin\frac{\pi}{k}}\log\frac{1}{\beta}+ c_k \log\log\frac{1}{\beta}+O(1),
\end{equation}
where $c_k$ is constant depending on $k$ only. Besides, unless
otherwise stated, the letter $C$ will always denote various
generic constants which are independent of $\beta$, especially for
$\beta$ sufficiently small. The constant $\alpha \in (0,
\frac{1}{2})$ is a fixed small constant.

The rest of this paper is organized as follows: In Section~2, we
introduce useful properties about the spike solution $w$ and the
symmetric vortex solution $v_d$. In Section~3, we define the
approximate solutions of spike-vortex solutions and derive the
associated error estimates. In Section~4, we use Liapunov-Schmidt
reduction process to find spike-vortex solutions. Then we may
complete the proof of Theorem~\ref{thm1_1} and \ref{thm1_2} in
Section~5 and 6, respectively. Finally, we give the proof of
Theorem~\ref{thm1_3} in Section~7.

\vskip 0.5cm

\noindent {\bf Acknowledgments:} The research of the first author
is partially supported by a research Grant from NSC of Taiwan.  The
research of the second author is partially supported by an
Earmarked Grant from RGC of Hong Kong.

\section{Properties of Spikes and Vortex}
\ \ \ \ We recall some properties of $w$ and
$S_d(r)e^{\sqrt{-1}\,d\theta}$, where $|d|=1$ or $|d|>1$ and
$2(d-1)\not\equiv 0 \mod k$. Let
\begin{equation}\label{id2_1}
\begin{cases}
L_1[\phi]=\triangle\phi-\phi+3w^2\phi,\\
L_2[\psi]=\triangle\psi+\psi-S_d^2\psi-2Re(S_d(r)e^{-\sqrt{-1}\,d\theta}
\psi)S_d(r)e^{\sqrt{-1}\,d\theta}\,,
\end{cases}
\end{equation} for $\phi$ is a real-valued function and $\psi$ is
a complex-valued function. For convenience, we may define the
conjugate operator of $L_2$ by
\begin{equation}
\hat{L}_2:= e^{-\sqrt{-1} d \theta } L_2 e^{\sqrt{-1} d \theta}
\end{equation}
Then by simple computations, it is easy to check that
\begin{equation}
\hat{L}_2 [ \psi_1 +\sqrt{-1} \psi_2] =  \hat{L}_{2,1} [\psi_1,
\psi_2] + \sqrt{-1}\,\hat{L}_{2,2} [\psi_1, \psi_2]\,,
\end{equation} for $\psi_1$ and $\psi_2$ are real-valued functions,
where
\[ \hat{L}_{2,1}[ \psi_1, \psi_2]= \Delta \psi_1 +(1-3S_d^2) \psi_1 -\frac{d^2}{r^2} \psi_1 -\frac{2 d}{r^2} \partial_\theta \psi_2\,,\]
\[ \ \hat{L}_{2,2}[ \psi_1, \psi_2]= \Delta \psi_2 +(1-S_d^2) \psi_2 -\frac{d^2}{r^2} \psi_2 +\frac{2 d}{r^2} \partial_\theta \psi_1\,.
\]
Set a function space
\begin{align}\label{id2_2}
\Sigma=\left\{
\left(
\begin{array}{cc}
\phi \\
\psi
\end{array} \right)=\left(
\begin{array}{cc}
\phi(z) \\
\psi(z)
\end{array} \right)
\in\mathbb{R}\times\mathbb{C}\Bigg|
\begin{array}{lll} &\phi(ze^{\sqrt{-1}\,\frac{2\pi}{k}})= \phi(z),
&\phi(\bar{z})=\phi(z),\\
&\psi(ze^{\sqrt{-1}\,\frac{2\pi}{k}})=e^{\sqrt{-1}\,\frac{2\pi}{k}}\psi(z),
&\psi(\bar{z})= \psi(z)^*\end{array} \right\}\,,
\end{align} where $k\geq 2$ is an integer. Hereafter,
both the over-bar and asterisk denote complex conjugate. We remark
that the equation (\ref{id1_5}) is invariant under the following
two maps
\begin{equation}\label{id2_3}
\begin{cases}
T_1
\left(
\begin{array}{cc}
\phi \\
\psi
\end{array} \right)
=
\left(
\begin{array}{cc}
\phi(ze^{\sqrt{-1}\,\frac{2\pi}{k}}) \\
e^{-\sqrt{-1}\,\frac{2\pi}{k}}\psi(ze^{\sqrt{-1}\,\frac{2\pi}{k}})
\end{array} \right),\\ \\
T_2
\left(
\begin{array}{cc}
\phi \\
\psi
\end{array} \right)
=
\left(
\begin{array}{cc}
\phi(\bar{z}) \\
\psi(\bar{z})^*
\end{array} \right)\,.
\end{cases}
\end{equation}
Therefore, we may look for solutions of (\ref{id1_5}) in the space
$\sum$. We first have
\begin{lem}\label{lem2_1}
\hspace*{0cm}
\begin{itemize}
\item[(1)]~~Suppose $L_1[\phi]=0$, $\phi\in H^2(\mathbb{R}^2)$ and
$\phi (\bar{z})=\phi(z)$. Then $\phi(z)=c\frac{\partial
w}{\partial z_1}(z)$, where $z=z_1+\sqrt{-1}z_2$,
$z_j\in\mathbb{R}$ and $c$ is a constant. \item[(2)]~~Suppose
\begin{equation}\label{id2_4} L_2[\psi]=0,\ |\psi|\leq C,\
\psi(z)^*=\psi(\bar{z}),\:\hbox{ and }\:
\psi(ze^{\sqrt{-1}\frac{2\pi}{k}})=e^{\sqrt{-1}\frac{2\pi}{k}}\psi(z)\,,
\end{equation} where $C$ is a positive constant.
Then $\psi\equiv 0$, provided that $|d|=1$ or $|d|>1$ and
$2(d-1)\not\equiv 0 \mod k$.
\end{itemize}
\end{lem}
\begin{proof}[\bf \large P{\small ROOF}]
(1) is easy to show. See Appendix C of \cite{nt2}.  We only need
to show (2). We firstly state the proof for the case when $|d|=1$.
For simplicity, we may assume $d=1$. From~\cite{PR}(Theorem~3.2),
we know that
\begin{align}\label{psi012}
\psi=c_0\big(\underbrace{\sqrt{-1}S(r)e^{\sqrt{-1}\theta}}_{\parallel\atop\displaystyle\psi_{0,0}}\big)+\sum_{j=1}^2c_j\underbrace{\frac{\partial}{\partial
z_j}\big(S(r)e^{\sqrt{-1}\theta}\big)}_{\parallel\atop\displaystyle\psi_{0,j}}\,,
\end{align} where $c_j$'s are constants.
It is easy to calculate that
$\left(\sqrt{-1}S(r)e^{\sqrt{-1}\theta}\right)^*=-\sqrt{-1}S(r)e^{-\sqrt{-1}
\theta}$,
$\psi_{0,1}=\big(\frac{S(r)}{r}\big)'\frac{x^2}{r}+\frac{S(r)}{r}+\sqrt{-1}\,(\frac{S(r)}{r})'xy$,
and
$\psi_{0,2}=\big(\frac{S(r)}{r}\big)'\frac{xy}{r}+\sqrt{-1}\big[(\frac{S(r)}{r})'y^2+\frac{S(r)}{r}\big]$,
where $x=z_1$ and $y=z_2$. Consequently,
$\left(\psi_{0,0}(z)\right)^*=-\psi_{0,0}(\bar{z})$,
$(\psi_{0,1}(z))^*=\psi_{0,1}(\bar{z})$ and
$(\psi_{0,2}(z))^*\neq\psi_{0,2}(\bar{z})$. Moreover, due to
$\psi(z)^*=\psi(\bar{z})$, we have $c_0=c_2=0$. Hence we only have
\begin{equation}\notag
\psi(z)=c_1\psi_{0,1}(z)\,.
\end{equation}
However, it is obvious that $\psi_{0,1}$ doesn't satisfy
$\psi_{0,1}\big(ze^{\sqrt{-1}\frac{2\pi}{k}}\big)=e^{\sqrt{-1}\frac{2\pi}{k}}\psi_{0,1}(z)$.
Thus $c_1=0$ and $\psi\equiv0$.

Now we give the proof for the case that $|d|>1$. From~\cite{L},
the solution $\psi$ locally may become a linear combination of
$\psi_{d,0}(z)= \sqrt{-1}\,h(r)\,e^{\sqrt{-1}\,d\theta}$ and the
following forms \beq\label{fm1} \psi_{d,m}(z)=
a(r)\,e^{\sqrt{-1}\,(d-m)\theta}+
b(r)\,e^{\sqrt{-1}\,(d+m)\theta}\,, \eeq for $m\in \mathbb{N}$,
where $z=r\,e^{\sqrt{-1}\,\theta}$ and $(r,\theta)$ is the polar
coordinate. Here $h$, $a$ and $b$ are real-valued functions.
Actually, these forms are invariant to the operator $L_2$ so one
may decompose the function space $\Sigma$ as invariant subspaces
having the forms like $\psi_{d,0}$ and $\psi_{d,m}$'s. Then the
condition $\psi(z)^*=\psi(\overline{z})$ may imply
$\psi_{d,0}(z)^*=\psi_{d,0}(\overline{z})$. However, since
$\psi_{d,0}(z)= \sqrt{-1}\,h(r)\,e^{\sqrt{-1}\,d\theta}$, then
$\psi_{d,0}(z)^*=-\psi_{d,0}(\overline{z})$ which gives
$\psi_{d,0}(z)\equiv 0$. Besides, the condition
$\psi(ze^{\sqrt{-1}\frac{2\pi}{k}})=e^{\sqrt{-1}\frac{2\pi}{k}}\psi(z)$
may give
$\psi_{d,m}(ze^{\sqrt{-1}\frac{2\pi}{k}})=e^{\sqrt{-1}\frac{2\pi}{k}}\psi_{d,m}(z)$.
Consequently, \beq\label{fm2}
a(r)\,e^{\sqrt{-1}\,(d-m)\theta}\left(e^{\sqrt{-1}\,(d-m)\frac{2\pi}{k}}
-e^{\sqrt{-1}\,\frac{2\pi}{k}}\right)
+b(r)\,e^{\sqrt{-1}\,(d+m)\theta}\left(e^{\sqrt{-1}\,(d+m)\frac{2\pi}{k}}-e^{\sqrt{-1}\,\frac{2\pi}{k}}\right)
=0\,.\eeq Hence $a\equiv 0$ or $b\equiv 0$ if
$e^{\sqrt{-1}\,(d-m-1)\frac{2\pi}{k}}\neq 1$ or
$e^{\sqrt{-1}\,(d+m-1)\frac{2\pi}{k}}\neq 1$. Due to
$L_2[\psi_{d,m}]=0$, $a(r)$ and $b(r)$ satisfy \beq\label{fm1.1}
\left\{\begin{array}{lll} &
a^{\prime\prime}+\frac{1}{r}a^\prime-\frac{(d-m)^2}{r^2}a+\left(1-S_d^2\right)a-(a+b)S_d^2
&=0\,, \\ &
b^{\prime\prime}+\frac{1}{r}b^\prime-\frac{(d+m)^2}{r^2}b
+\left(1-S_d^2\right)b-(a+b)S_d^2 &=0\,.\end{array}\right.\eeq
This implies $a\equiv b\equiv 0$ if $a\equiv 0$ or $b\equiv 0$.
Thus $a\equiv b\equiv 0$ if
$e^{\sqrt{-1}\,(d-m-1)\frac{2\pi}{k}}\neq 1$ or
$e^{\sqrt{-1}\,(d+m-1)\frac{2\pi}{k}}\neq 1$. It is trivial that
$e^{\sqrt{-1}\,(d-m-1)\frac{2\pi}{k}}\neq 1$ or
$e^{\sqrt{-1}\,(d+m-1)\frac{2\pi}{k}}\neq 1$ if $2(d-1)\not\equiv
0\mod k$. Therefore $a\equiv b\equiv 0$ i.e. $\psi_{d,m}\equiv 0$
if $2(d-1)\not\equiv 0\mod k$, and we may complete the proof of
Lemma~\ref{lem2_1}.
\end{proof}

To study the properties of $L_2$ (or $\hat{L}_2$), we introduce
some Sobolev spaces. Let $\alpha\in(0,\frac{1}{2})$. We introduce
Hilbert spaces $X_\alpha$ and $Y_\alpha$ as follows:
\begin{equation}\notag
X_\alpha=\Bigg\{u= u_1 +\sqrt{-1}\, u_2\in
L^2_{loc}(\mathbb{R}^2;\mathbb{C})\Bigg|\int_{\mathbb{R}^2}(1+
|x|^{2+\alpha})(|u_1|+|u_2|)^2<+\infty\Bigg\}
\end{equation}
equipped with inner product
$\displaystyle(u,v)_{X_\alpha}=\int_{\mathbb{R}^2}(1+|x|^{2+\alpha})
(u_1 v_1+ u_2 v_2)dx$, and
\begin{equation}\notag
Y_\alpha=\Bigg\{v= v_1 +\sqrt{-1}\,v_2\in
W_{loc}^{2,2}(\mathbb{R}^2;\mathbb{C}) \Bigg|\int_{\mathbb{R}^2
}|\triangle
v|^2(1+|x|^{2+\alpha})dx+\int_{\mathbb{R}^2}\frac{|v|^2}{1+|x|^{2+\alpha}}dx<+\infty\Bigg\}
\end{equation}
equipped with inner product $\displaystyle(u,v)_{Y_\alpha}=(\triangle u,
\triangle v)_{X_\alpha}+\int_{\mathbb{R}^2}\frac{u\cdot v}{1+|x|^{2+\alpha}}dx$,
respectively. Thanks to the inequality
\[ \int_{\R^2} |h| \leq \left(\int_{\R^2} \frac{1}{ (1+|z|)^{2+\alpha}}\right)^{\frac{1}{2}} \|
h\|_{X_\alpha}\,.\] Besides, we see that $ X_\alpha $ has a
compact embedding to $L^1 (\R^2)$. Originally, these spaces are
real-valued function spaces introduced in Chae and Imanuvilor
(cf.~\cite{DI}). Here we generalize $X_\alpha,\ Y_\alpha$ as
complex-valued function spaces, and regard the operator $L_2$ from
the space $Y_\alpha$ to the space $X_\alpha$. We list some
properties of $X_\alpha$ and $Y_\alpha$, whose proofs are exactly
the same as in~\cite{DI}.
 \begin{lem}\label{lem2_2}
 \hspace*{0cm}
\begin{itemize}
\item[(1)] Let $v\in Y_\alpha$ be a harmonic function.
 Then $v=$constant.
 \item[(2)]$\forall v\in Y_\alpha$, we have $|v(z)|\leq C_1||v||_{Y_\alpha}(\ln(1+|z|)+1)$, $\forall z\in\mathbb{R}^2$.
\item[(3)] The image of $L_2$ (or $\hat{L}_2$) is closed in
$X_\alpha \cap \Sigma_0$, where $\Sigma_0\equiv
\{\psi=\psi(z)\in\mathbb{C}: (0,\psi)\in \Sigma\}\,.$
\end{itemize}
\end{lem}

Now we study the invertibility of $L_2$ (or $\hat{L}_2$) on the
space $Y_\alpha\cap\Sigma_0$.
\begin{lem}\label{lem2_3}
For $\alpha\in(0,\frac{1}{2})$. Then operator $L_2$ (or
$\hat{L}_2$) from the space $Y_\alpha\cap\Sigma_0$ onto the space
$X_\alpha\cap\Sigma_0$ is invertible. Furthermore,
\begin{equation}
\| \psi\|_{Y_\alpha} \leq C \| L_2[\psi]\|_{X_\alpha}\,, \quad \|
\psi\|_{Y_\alpha} \leq C \| \hat{L}_2[\psi]\|_{X_\alpha}
\end{equation}
\end{lem}
\begin{proof}[\bf \large P{\small ROOF}]
%Actually, we only need to show that $$\displaystyle
%Im(\hat{L}_2)=\Bigg\{f\in X_\alpha\Bigg| \mbox{Re}
%\left(\int_{\mathbb{R}^2}f \psi_{0,j} e^{-\sqrt{-1} d
%\theta}\right)=0, j=0, 1, 2\Bigg\}\,.$$ (Since $\psi_{0,j} e^{
%-\sqrt{-1} d \theta} \not \in\Sigma_0,\ j=0,1,2$,  this then
%implies that the image of $\hat{L}_2$ from $Y_\alpha \cap
%\Sigma_0$ is $ X_\alpha \cap \Sigma_0$, which proves the
%ontoness.)
It suffices to consider the invertibility of $\hat{L}_2$. Then the
invertibility of $L_2$ follows from a trivial transformation. To
show the invertibility of $\hat{L}_2$, we claim that
$(Im(\hat{L}_2))^\perp= \{0\}\,.$ Suppose
$\psi\in(Im(\hat{L}_2))^\perp$. Then it is easy to check that
$\hat{L}_2 [\psi]=0, \psi \in X_\alpha\cap\Sigma_0$. By Lemma
\ref{lem2_1}, we just need to show that $\psi$ is bounded. To this
end, we note that since $ \psi\in \Sigma_0$, we have $\psi (0)=
e^{ \sqrt{-1} \frac{2\pi}{k}} \psi (0)$ and $k\geq 2$ which may
give $\psi(0) =0\,.$ Let $\phi= \psi\,e^{\sqrt{-1}\,d\theta}$.
Then $\phi\in X_\alpha$ and $L_2 [\phi]=0$. Hence $\phi\in
Y_\alpha$ and
$$\triangle \phi(z)= h(z)/\left(1+|z|^{2+\alpha}\right)^{1/2}\,,
\forall z\in\mathbb{R}^2\,,$$ where $h\in L^2(\mathbb{R}^2)$.
Consequently, by Lemma~\ref{lem2_2} and Riesz representation
formula, we may obtain  $$|\phi(z)| \leq C \log (1+|z|)\,,\quad
\forall z\in\mathbb{R}^2\,,$$ and $$|\partial_\theta\,\phi(z)|\leq
C \log (1+|z|)\,,\quad \forall z\in\mathbb{R}^2\,.$$ Thus
\beq\label{r.1} |\psi(z)|\,,\: |\partial_\theta \psi(z) | \leq C
\log (1+|z|)\,,\quad \forall z\in\mathbb{R}^2\,.\eeq

The equation $\hat{L}_2 [\psi]=0, \psi=\psi_1+\sqrt{-1}\psi_2$ can
be written as \begin{eqnarray}\label{1} \Delta \psi_1 +(1-3S_d^2)
\psi_1 -\frac{d^2}{r^2} \psi_1 -\frac{2 d}{r^2} \partial_\theta
\psi_2 &=& 0\,, \\ \label{2}\Delta \psi_2 +(1-S_d^2) \psi_2
-\frac{d^2}{r^2} \psi_2 +\frac{2 d}{r^2} \partial_\theta \psi_1
&=& 0\,.
\end{eqnarray}
From (\ref{r.1}) and (\ref{1}), $\psi_1$ satisfies
$$
 |\Delta \psi_1 + (1-3S_d^2)\psi_1 - \frac{d^2}{r^2}\psi_1|\leq
 \frac{C}{r^2}\,\log r\quad\hbox{ for }\:r=|z|>1 \,.
$$
Hence by comparison principle, \beq\label{2.5}|\psi_1(z)| \leq
\frac{C}{1+|z|}\,.\eeq Here we have used the fact that $S_d(r)\sim
1$ as $r\to\infty$. Similarly, we may have \beq\label{3}
|\partial_\theta \psi_1(z)|\leq \frac{C}{1+|z|}\,.\eeq

It remains to show that $\psi_2 $ is bounded. Now we can use
(\ref{2}) and (\ref{3}) to get $$ |\Delta \psi_2 + (1-S_d^2)\psi_2
- \frac{d^2}{r^2}\psi_2|\leq
 \frac{C}{r^2(1+r)}\quad\hbox{ for }\:r=|z|>1 \,.
$$ In fact, $\psi_2$
satisfies
\[ -\Delta \psi_2 = f(z), \ \mbox{where} \ f(z)= (1-S_d^2-\frac{d^2}{r^2}) \psi_2 +\frac{2 d}{r^2} \partial_\theta \psi_1\]
Since $\psi \in \Sigma_0$, we deduce that $ f(z) \in X_\alpha \cap
\Sigma_0$ and $ \int_{\R^2} f(z)=0$. Since $\psi_2 (0)=0$, we have
\begin{equation}
\label{psi2}
\psi_2 (z)=\psi_2 (0)+\frac{1}{2\pi} \int_{\R^2} \log \frac{|\tau|}{|z-\tau|}  f(\tau) d \tau = \frac{1}{2\pi} \int_{\R^2} \log \frac{|\tau| |z| }{|z-\tau|}  f(\tau) d \tau
\end{equation}
from which we obtain that \beq\label{4} |\psi_2(z)| \leq C \| f
\|_{X_\alpha} <+\infty\quad\hbox{ for }\:z\in\mathbb{R}^2\,. \eeq
Therefore by (\ref{2.5}) and (\ref{4}), we may complete the proof
of Lemma~\ref{lem2_3}.
\end{proof}

To study the vortex solutions,  we perform the following key transformation:
\begin{equation}
v= v_d (z)  e^{ \sqrt{-1} \psi (z)}\,,\quad v_d(z)=
S_d(r)\,e^{\sqrt{-1}\,d\theta}\,,
\end{equation} where $\psi=\psi_1
+\sqrt{-1} \psi_2$ and $\psi_j\in\mathbb{R}\,, j=1,2\,.$  (Here we have {\it assumed} that the vortex occurs at the center.) Now we define
\[ S_1[u, \psi]:= \Delta u-u +u^3 + \beta u  S_d^2  e^{- 2 \psi_2}\,, \]
\[ S_2 [u, \psi]:= \Delta \psi + \frac{2 \nabla v_d}{v_d}\cdot\nabla \psi
- \sqrt{-1} S_d^2 (1-e^{-2\psi_2}) + \sqrt{-1} |\nabla \psi|^2-
\sqrt{-1} \beta u^2\,.  \] We also write
\begin{equation}
 S_2 [u, \psi]= S_2[u, 0]+ \tilde{L}_2 [\psi] + N [u,\psi]
\end{equation}
where
\begin{eqnarray}
\label{newL}
 \tilde{L}_2 [\psi] &=& \Delta \psi +\frac{2 \nabla v_d}{v_d} \nabla \psi- 2 \sqrt{-1} S_d^2
 \psi_2 \\ \nonumber
 &=& \left[\Delta \psi_1 + 2 \left( \frac{ S_d^{'} (|z|)}{S_d
(|z|)} \frac{z}{|z|}\right) \nabla \psi_1
-d\nabla\theta\cdot\nabla\psi_2\right] \\ \label{2.22} & &~~+
\sqrt{-1} \left[\Delta \psi_2 -2 S_d^2 \psi_2 + 2
\left(\frac{S_d^{'}}{S_d} \frac{ z}{|z|}\right) \nabla \psi_2
+d\nabla\theta\cdot\nabla\psi_1 \right]\,, \end{eqnarray}
\[ N [u, \psi]= \sqrt{-1} |\nabla \psi|^2 -\sqrt{-1} S_d^2 (1- e^{-2 \psi_2} -2\psi_2)\,,\]
where $\nabla\theta=\frac{1}{r}(-\sin\theta, \cos\theta)$. Then it
is easy to see that solving (\ref{id1_5}) is equivalent to solving
\begin{equation}
\label{neweqn} S_1 [u, \psi]=0,\quad S_2 [u, \psi]=0\,.
\end{equation}

Let $ \psi=\psi_1+ \sqrt{-1}  \psi_2, h= h_1 + \sqrt{-1} h_2$. We
may define two norms as follows:
\begin{equation}
\label{norm} \| \psi \|_{*}= \sup_{z \in \R^2} [  | \psi_1 |+
(1+|z|) |\nabla \psi_1|+ (1+|z|)^{1+\alpha} | \psi_2 |+
(1+|z|)^{1+\alpha} |\nabla \psi_2| ]\,,
\end{equation} and
\[\| h \|_{**}= \sup_{z \in \R^2} [(1+|z|)^{2+\alpha} (|h_1|+| h_2
|)]\,.
\]

Then we may show the following key lemma
\begin{lem}
\label{lem2_4} For any $ h \in X_\alpha \cap \Sigma_0$ with $ \|
h\|_{**} <+\infty$,  there exists a unique $\tilde{\psi} \in
Y_\alpha \cap \Sigma_0$ such that
\begin{equation}
\label{lh}
\tilde{L}_2 [\psi]= h
\end{equation}
Furthermore, we have
\begin{equation}
\label{phi2}
\|\psi\|_{*} \leq C \| h \|_{**}
\end{equation}
\end{lem}

\noindent {\bf Proof:} Let $ \hat{h}= -\sqrt{-1}\,v_d h$ and
$\hat{ \phi} = -\sqrt{-1}\, v_d \psi\,.$ Then (\ref{lh}) is
equivalent to
\begin{equation}
\label{phi1}
L_2 [ \hat{\phi} ]= \hat{h}
\end{equation}
where $\hat{h}$ satisfies $ \hat{h} \in \Sigma_0$ and
\begin{equation}
\label{hin} \|\hat{ h} \|_{**}\leq  \|h\|_{**} <+\infty\,.
\end{equation}
Consequently, (\ref{hin}) implies that $ \int_{\R^2} |\hat{h}|^2
(1+ |z|^{2+\alpha}) <+\infty$ and hence  $\hat{ h }\in X_\alpha
\cap \Sigma_0$. By Lemma \ref{lem2_3}, (\ref{phi1}) has a unique
solution $ \phi$. Furthermore, as for (\ref{psi2}), we obtain
\begin{equation}
\label{phi} \hat{\phi}(z)= \frac{1}{2\pi} \int_{\R^2} \log \frac{
|\tau| |z|}{ |z-\tau|} \hat{f} (\tau) d \tau
\end{equation}
where $ |\hat{f}(z)|\,( 1+|z|)^{2+\alpha} \leq C \| h\|_{**} $.
From (\ref{phi}), we may deduce that $ |\hat{\phi} | + (1+|z|)
|\nabla \hat{\phi} | \leq C \| h\|_{**}$ which implies that $
|\psi| + (1+|z|) |\nabla \psi| \leq C \| h\|_{**}$ i.e.~$ |\psi_j
| + (1+|z|) |\nabla \psi_j| \leq C \| h\|_{**}\,, j=1,2\,.$

To obtain better estimates for $\psi_2$, we may use the equation
for $\psi_2$ and (\ref{2.22}). Then
\[
\Delta \psi_2 -2 S_d^2 \psi_2 + O\left( \frac{1}{|z|} |\nabla
\psi|\right) = O\left(\|h\|_{**}\,( 1+|z|)^{-2-\alpha}\right)\,,
\]
which gives
$$
-\frac{C\,\|h\|_{**}}{ (1+|z|)^2} \leq - \Delta \psi_2 +2 S_d^2
\psi_2 \leq \frac{C\,\|h\|_{**}}{ (1+|z|)^2}\,.
$$
Hence we may use a barrier and elliptic estimates to get
$$
 |\psi_2| +|\nabla \psi_2| \leq \frac{C\,\|h\|_{**}}{ (1+|z|)^{1+\alpha}}\,.
$$ Here we have used the fact that $S_d(r)\sim 1$ as $r\to \infty$.
Therefore we may complete the proof of Lemma~\ref{lem2_4}.
\qed

\section{Approximate Solutions and Error Estimates}
\ \ \ \ In this section, we introduce some approximate solutions
and derive some useful estimates. Let
\begin{eqnarray*}
S[u,\psi]=&
\left(
\begin{array}{cc}
S_1[u, \psi] \\
S_2[u, \psi]
\end{array} \right).
\end{eqnarray*}
Let
\begin{align}
\xi_j=le^{\sqrt{-1}\frac{2\pi}{k}(j-1)}&,\ j=1,...,k,\ w_j(z):=w(z-\xi_j
),\notag\\
&u_l(z)=\sum_{j=1}^kw_j(z),\notag
\end{align}
Then we have
\begin{lem}\label{lem3_1}
For $l$ large enough, we have
\begin{equation}\label{id3_1}
\big\|S_1[u_l, 0]\big\|_{L^2(\mathbb{R}^2)}+\big\|S_2[u_l,
0]\big\|_{**} \leq C\,\big(|\beta|
l^{2+\alpha}+e^{-2l\sin\frac{\pi}{k}}\big)\,.
\end{equation}
\end{lem}
\begin{proof}[\bf \large P{\small ROOF}]
It is easy to check that
\begin{align}
S_1[u_l, 0]=&\triangle u_l-u_l+u_l^3+\beta u_l |v_d|^2 \notag\\
=&\bigg(\sum_jw_j(z)\bigg)^3-\sum_jw_j^3+\beta\sum_j w_j S_d^2\notag\\
=&O\bigg(\sum_{i\neq j}w_i^2w_j+|\beta|\sum_j|w_j|\bigg).\notag
\end{align}
Due to
\begin{equation}
\int_{\R^2} w_i^4w_j^2\leq Ce^{-2|\xi_i-\xi_j|} \leq
Ce^{-4l\sin\frac{\pi}{k}}\,,\quad \forall i\neq j\,,
\end{equation}
we may obtain
\begin{equation}\label{id3_2}
\big\|S_1[u_l,
0]\big\|_{L^2}\leq\big[e^{-2l\sin\frac{\pi}{k}}+|\beta|\big]\,.
\end{equation} Here we have used the fact that $\xi_j$'s are
vertices of regular $k$-polygon with side length
$2l\sin\frac{\pi}{k}\,.$ On the other hand, we have
\[ |S_2 [u_l, 0]| = |\beta| u_l^2 \leq C \sum_{j=1}^k |\beta| e^{-2 |z-\xi_j|}
\]
 and so by (\ref{ldef})
\begin{equation}
\label{id3_3}
 \| S_2 [u_l, 0] \|_{**}  \leq C |\beta| \sup_{z \in \R^2}
 \left(\sum_{j=1}^k\,|z|^{2+\alpha}\,e^{-2 |z-\xi_j|}\right)
 \leq C | \beta| \sum_{j=1}^k |\xi_j|^{2+\alpha} \leq C |\beta|
 |l|^{2+\alpha}\,.
\end{equation}
Therefore by (\ref{id3_2}) and (\ref{id3_3}), we may obtain
(\ref{id3_1}) and complete the proof of Lemma~\ref{lem3_1}.
\end{proof}

\section{Liapunov-Schmidt reduction process}

\ \ \ \ Let $\widetilde{X}:=(L^2(\mathbb{R}^2)\times
X_\alpha)\cap\Sigma$, $\widetilde{Y}:=(H^2(\mathbb{R}^2)\times
Y_\alpha)\cap\Sigma$ and $L:= \left( {\begin{array}{*{20}c}
    \widetilde{L}_1 \\\widetilde{L}_2
\end{array}} \right):\widetilde{Y}\rightarrow\widetilde{X}\,,
$ where $\widetilde{L}_1[\phi]\equiv\triangle \phi -\phi
+3u_l^2\,\phi$ and $\ \widetilde{L}_2$ is defined by (\ref{newL}).
To solve (\ref{neweqn}), we first consider the following linear
problem: Given $f\in L^2(\R^2)\cap\Sigma_1$, find $(\phi,c)$ such
that
 \begin{equation}\label{id4_1}
\begin{cases}
\displaystyle\widetilde{L}_1\phi=f+c\frac{\partial u_l}{\partial l}\,,\quad
\phi\in H^2(\R^2)\cap\Sigma_1,\\
\displaystyle\int_{\R^2} \phi\frac{\partial u_l}{\partial l}=0\,,
\end{cases}
 \end{equation} where $\Sigma_1\equiv \{\phi=\phi(z)\in\mathbb{R}:
 (\phi, 0)\in\Sigma\}$.
 \begin{lem}\label{lem4_1}
 For each $f\in L^2(\R^2)\cap\Sigma_1$, there exists a unique pair $(\phi,c)$ satisfying (\ref{id4_1}) such that
 \begin{equation}\label{id4_2}
\big\|\phi\big\|_{H^2}\leq C\big\|f\big\|_{L^2}.
\end{equation}
\end{lem}
\begin{proof}[\bf \large P{\small ROOF}]
The existence and uniqueness may follow from Fredholm's
alternatives. Now we prove (\ref{id4_2}) by contradiction. Suppose
(\ref{id4_2}) is not true. Then there exist $\beta_n$, $l_n\in
\mathbb{R}$, $\phi_n\in H^2(\R^2)\cap\Sigma_1\,,$ $f_n\in
L^2(\R^2)\cap\Sigma_1$ and $c_n$ such that $\beta_n\rightarrow0,\
l_n\rightarrow+\infty\,,$
\begin{equation}\label{id4_3}
\big\|f_n\big\|_{L^2}\rightarrow0,\quad
\big\|\phi_n\big\|_{H^2}=1\,,
\end{equation} as $n\to\infty$,
\beq\label{4.1.1}
\displaystyle\widetilde{L}_1\phi_n=f_n+c_n\frac{\partial
u_l}{\partial l}\,, \eeq and \beq\label{4.1.2}
\displaystyle\int_{\R^2} \phi_n\frac{\partial u_l}{\partial
l}=0\,. \eeq Let $\widetilde{\phi}_n(z):=\phi_n(z+\xi_1)$. Then by
(\ref{4.1.1}), $\widetilde{\phi}_n$ satisfies
\beq\label{id4_5}\begin{array}{lll}
&\triangle\widetilde{\phi}_n(z)-\widetilde{\phi}_n(z)
+3w^2(z)\widetilde{\phi}_n(z) \\
&+3\left[2\ds\sum_{j=2}^k\,w(z)w(z+\xi_1-\xi_j)
+\left(\ds\sum_{j=2}^k\,w(z+\xi_1-\xi_j)\right)^2\right]\widetilde{\phi}_n(z) \\
&=f_n(z+\xi_1)+c_n\Bigg[-\frac{\partial w}{\partial
z_1}(z)+\ds\sum_{j=2}^k\frac{\partial w_j}{\partial
l}(z+\xi_1)\Bigg]\,.\end{array} \eeq We may multiply (\ref{4.1.1})
by $\displaystyle\frac{\partial u_l}{\partial l}$ and integrate
over the whole space $\mathbb{R}^2$. Then by (\ref{id4_3}), it is
obvious that
\begin{align}
c_n\int_{\mathbb{R}^2}\bigg(\frac{\partial u_l}{\partial l}\bigg)^2
&=\int_{\mathbb{R}^2}\big(\triangle\phi_n-\phi_n+3u_l^2\phi_n\big)
\frac{\partial u_l}{\partial l}-\int_{\mathbb{R}^2}f_n\frac{\partial u_l}{\partial l}\notag\\
&=\int_{\mathbb{R}^2}\bigg(\triangle\frac{\partial u_l}{\partial l}
-\frac{\partial u_l}{\partial l}+3u_l^2\frac{\partial u_l}{\partial l}\bigg)\phi_n
-\int_{\mathbb{R}^2}f_n\frac{\partial u_l}{\partial l}\notag\\
&\stackrel{n\rightarrow+\infty}{-\!\!\!-\!\!\!-\!\!\!-\!\!\!\longrightarrow}0.\notag
\end{align}
So $c_n\rightarrow 0$ as $n\rightarrow+\infty$. Hence by
(\ref{id4_5}), as $n\rightarrow\infty$,
$\widetilde{\phi}_n\rightarrow\widetilde{\phi}_0$ which satisfies
$\triangle\widetilde{\phi}_0-\widetilde{\phi}_0+3w^2\widetilde{\phi}_0=0$
in $\mathbb{R}^2$. Thus by~\cite{nt2},
$\displaystyle\widetilde{\phi}_0=\sum_{j=1}^2\,a_j\frac{\partial
w}{\partial z_j}$ for some constants $a_1$ and $a_2$. Moreover,
due to $\widetilde{\phi}_n(\bar{z})=\widetilde{\phi}_n(z)$, we
have $\widetilde{\phi}_0(\bar{z})=\widetilde{\phi}_0(z)$.
Consequently, $a_2=0$ and
$\displaystyle\widetilde{\phi}_0=a_1\frac{\partial w}{\partial
z_1}$. On the other hand, by (\ref{4.1.2}),
\begin{align}
0&=\int_{\mathbb{R}^2}\phi_n\frac{\partial u_l}{\partial l}\notag\\
&=\int_{\mathbb{R}^2}\widetilde{\phi}_n(z)\left[-\frac{\partial w}{\partial z_1}(z)
+ \ds\sum_{j=2}^k\,\frac{\partial w_j}{\partial l}(z+\xi_1)\right]\notag\\
&\stackrel{n\rightarrow+\infty}{-\!\!\!-\!\!\!-\!\!\!-\!\!\!\!\longrightarrow}
-a_1\int_{\mathbb{R}^2}\Bigg(\frac{\partial w}{\partial
z_1}\Bigg)^2\,.\notag
\end{align}
Therefore $a_1=0$, $\widetilde{\phi}_0=0$, and
$\widetilde{\phi}_n\rightarrow 0$ in $L^2_{loc}(\mathbb{R}^2)$.
Then we have $3u_l^2\widetilde{\phi}_n\rightarrow0$ in $L^2$ and
hence $\|\phi_n\|_{H^2}\rightarrow0$ which may give a
contradiction and complete the proof.
\end{proof}
From Lemma~\ref{lem2_4} and \ref{lem4_1}, we may obtain that
\begin{lem}\label{lem4_2}
For $\left( {\begin{array}{*{20}c}
    f_1 \\f_2
\end{array}} \right)
\in\widetilde{X},\ \mbox{there exists a unique} \
\bigg(\left( {\begin{array}{*{20}c}
    \phi \\ \psi
\end{array}} \right),c\bigg)
\in\widetilde{Y}\times\mathbb{R}
$, such that
\begin{align}\label{id4_6}
\widetilde{L}
\left( {\begin{array}{*{20}c}
    \phi \\ \psi
\end{array}} \right)
=
\left( {\begin{array}{*{20}c}
    f_1 \\f_2
\end{array}} \right)
+c
\left( {\begin{array}{*{20}c}
    \displaystyle\frac{\partial u_l}{\partial l} \\0
\end{array}} \right)
\,,\quad \int_{\mathbb{R}^2}\phi\displaystyle\frac{\partial
u_l}{\partial l}=0.
\end{align}
Moreover, we have $\|\phi\|_{H^2}+\|\psi\|_{*}\leq
C\,(\|f_1\|_{L^2}+\|f_2\|_{**})$.
\end{lem}
We may denote $A$ as the inverse operator for Lemma~\ref{lem4_2},
i.e. $A \left( {\begin{array}{*{20}c}
    f_1 \\f_2
\end{array}} \right)
=
\left( {\begin{array}{*{20}c}
    \phi \\ \psi
\end{array}} \right).
$ Finally, we have
\begin{lem}\label{lem4_3}
For $l$ large satisfying (\ref{ldef}), there exists a unique
$
\left( {\begin{array}{*{20}c}
    \phi_l \\ \psi_l
\end{array}} \right)
\in \widetilde{Y} $ such that
\begin{equation}\label{id4_7}
S[u_l+\phi_l,\psi_l]=c(l)
\left( {\begin{array}{*{20}c}
    \displaystyle\frac{\partial u_l}{\partial l} \\ 0
\end{array}} \right)
\,,\quad \int_{\mathbb{R}^2}\phi_l\displaystyle\frac{\partial
u_l}{\partial l}=0.
\end{equation}
Furthermore,
\begin{equation}\notag
\|\phi_l\|_{H^2}+\|\psi_l\|_{*}\leq C\,\Big(
e^{-2l\sin\frac{\pi}{k}}+|\beta|l^{2+\alpha}\Big).
\end{equation}
\end{lem}
\begin{proof}[\bf \large P{\small ROOF}]
This may follow from standard contraction mapping principle. We
choose $(\phi, \psi)\in B$, where $$
B=\left\{(\phi,\psi)\in\widetilde{Y}:
\|\phi\|_{H^2}+\|\psi\|_{*}\leq
C\,\big(e^{-2l\sin\frac{\pi}{k}}+|\beta|
l^{2+\alpha}\big)\right\}\,,
$$ and then expand
\begin{align}
S_1[u_l+\phi,\psi]=S_1[u_l, 0]+\tilde{L}_1[\phi]+N_1[\phi,\psi],\notag\\
S_2[u_l+\phi, \psi]=S_2[u_l,0]+\tilde{L}_2[\psi]+N_2[\phi,\psi],\notag
\end{align}
where
\begin{equation}\label{4.9}
%\big|N_1[\phi,\psi]\big|\leq C\,\big(\beta|u_l||\phi|+ u_l
%|\phi|^2+ \beta u_l^2|\psi_2|\big)
N_1[\phi,\psi]= 3u_l\phi^2 +\phi^3
+\beta\left[\phi\,S_d^2\,e^{-2\psi_2} +
u_l\,S_d^2\left(e^{-2\psi_2}-1\right)\right]\,,
\end{equation}
and
\begin{equation}\label{4.10}
N_2 [\phi, \psi]= -\sqrt{-1}\,\beta (2 u_l \phi +\phi^2) +
\sqrt{-1}\,|\nabla \psi|^2 - \sqrt{-1} S_d^2 (1-e^{-2\psi_2}
-2\psi_2)\,.
\end{equation}
By~(\ref{ldef}) and $(\phi,\psi)\in B$, we have $$
\|\phi\|_{H^2}+\|\psi\|_{*}\leq
|\beta|^{1-\sigma}\,,\quad\sigma\in
\left(0\,,\frac{1}{100}\right)\,,
$$ as $\beta>0$ sufficiently small i.e. $l$ sufficiently large.
Moreover, by (\ref{4.9}), (\ref{4.10}) and the norms defined at
(\ref{norm}), we see that
\begin{equation}\label{4.11}
\| N_1\|_{L^2(\R^2)} \leq C |\beta| \|\phi\|_{L^2
}+\|\phi\|_{L^2}^2 + |\beta|\,,
\end{equation}
and
\begin{equation}\label{4.12}
\|N_2 \|_{**} \leq C \| \psi\|_{*}^2 + C |\beta|\,.
\end{equation}

Now we can write (\ref{id4_7}) as
\begin{equation}\notag
\left( {\begin{array}{*{20}c}
    \phi \\ \psi
\end{array}} \right)
=A
\left( {\begin{array}{*{20}c}
    S_1[u_l, 0]+N_1[\phi,\psi] \\ S_2[u_l, 0]+N_2[\phi,\psi]
\end{array}} \right).
\end{equation}
Then as for the proof of Proposition~1 in \cite{LW2}, we may use a
contraction mapping argument to obtain the desired result. Here we
also need to use (\ref{4.11}), (\ref{4.12}), Lemma~\ref{lem3_1},
\ref{lem4_1} and \ref{lem4_2}.
\end{proof}
\section{Finding zero of $c(l)$}
\ \ \ \  To prove Theorem~\ref{thm1_1}, it is enough to find a
zero of $c(l)$ in~(\ref{id4_7}). We multiply the first equation of
(\ref{id4_7}) by $\displaystyle\frac{\partial u_l}{\partial l}$
and integrate it over the whole space $\mathbb{R}^2$. Then we
obtain
\begin{align}\label{id5_1}
c(l)\int_{\mathbb{R}^2}\bigg(\frac{\partial u_l}{\partial
l}\bigg)^2=&\int_{\mathbb{R}^2}\big[\triangle(u_l+\phi_l)-(u_l+\phi_l)+(u_l+\phi_l)^3\big]\frac{\partial
u_l}{\partial l} \\
&+\beta\int_{\mathbb{R}^2}(u_l+\phi_l)|v_d +\psi_l|^2\frac{\partial u_l}{\partial l}\notag\\
=&:I_1+I_2,
\end{align}
where
\begin{align}
I_2=&\beta\int_{\mathbb{R}^2}u_l|v_d|^2\frac{\partial u_l}{\partial l}+\beta\int_{\mathbb{R}^2}
O\left(|\phi_l|+|u_l+\phi_l||\psi_l|\right)\bigg|\frac{\partial u_l}{\partial l}\bigg|\notag\\
=&\beta\int_{\mathbb{R}^2}u_lS_d^2\frac{\partial u_l}{\partial l}+O\Big(|\beta|^2l^{2+\alpha}+|\beta|e^{-2l\sin\frac{\pi}{k}}\Big)\notag\\
=&\beta\int_{\mathbb{R}^2}u_l\Bigg(1-\frac{d^2}{2r^2}+O\bigg(\frac{1}{r^4}\bigg)\Bigg)^2\frac{\partial u_l}{\partial l}+O\Big(|\beta|^2l^{2+\alpha}+|\beta|e^{-2l\sin\frac{\pi}{k}}\Big)\notag\\
=&\beta\int_{\mathbb{R}^2}u_l\frac{\partial u_l}{\partial l}+O\Big(|\beta|l^{-4}+|\beta|^2l^{2+\alpha}+|\beta|e^{-2l\sin\frac{\pi}{k}}\Big)-\beta\int_{\mathbb{R}^2}u_l\frac{d^2}{r^2}\frac{\partial u_l}{\partial l}.\notag
\end{align}
Note that
\begin{align}\label{id5_2}
\int_{\mathbb{R}^2}u_l\frac{\partial u_l}{\partial l}=&\int_{\mathbb{R}^2}\sum_{j=1}^k w_j\sum_{i\neq j}\frac{\partial w_i}{\partial l}\notag\\
%=&\sum_j\int_{\mathbb{R}^2}w_j\frac{\partial u_l}{\partial l}+O\Big(e^{-2l\sin\frac{\pi}{k}}\Big)\notag\\
=&O\Big(e^{-2l\sin\frac{\pi}{k}}\Big),
\end{align}
\begin{align}\label{id5_3}
-\beta\int_{\mathbb{R}^2}u_l\frac{d^2}{r^2}\frac{\partial u_l}{\partial l}
=&\beta\left[\int_{\mathbb{R}^2}w\frac{d^2}{|z+\xi_1|^2}\frac{\partial w}{\partial z_1}+O\Big(e^{-2l\sin\frac{\pi}{k}}\Big)\right]\notag\\
=&-c_2\,\beta
l^{-3}\int_{\mathbb{R}^2}w^2+O\Big(|\beta|l^{-5}+|\beta|e^{-2l\sin\frac{\pi}{k}}\Big)\,,
\end{align} where $c_2$ is a positive constant independent of
$\beta$ and $l$. So
\begin{align}\label{id5_4}
I_2=-c_2\,\beta
l^{-3}\int_{\mathbb{R}^2}w^2+O\Big(|\beta|l^{-4}+|\beta|^2l^{2+\alpha}+|\beta|e^{-2l\sin\frac{\pi}{k}}\Big).
\end{align}
For $I_1$, we have
\begin{align}\label{id5_5}
I_1=&\int_{\mathbb{R}^2}\big(\triangle u_l-u_l+u_l^3\big)\frac{\partial u_l}{\partial l}
+\int_{\mathbb{R}^2}\big(\triangle \phi_l-\phi_l+3u_l^2\phi_l\big)\frac{\partial u_l}{\partial l}
+O\Big(|\beta|^2l^{2+\alpha}+|\beta|e^{-2l\sin\frac{\pi}{k}}\Big)\notag\\
=&\int_{\mathbb{R}^2}\bigg[\bigg(\sum_jw_j\bigg)^3-\sum_jw_j^3\bigg]\frac{\partial u_l}{\partial l}
+\int_{\mathbb{R}^2}\bigg[\triangle\frac{\partial u_l}{\partial l}-\frac{\partial u_l}{\partial l}+3u_l^2\frac{\partial u_l}{\partial l}\bigg]\phi_l
+O\Big(|\beta|^2l^{2+\alpha}+|\beta|e^{-2l\sin\frac{\pi}{k}}\Big)\notag\\
=&6k\int_{\mathbb{R}^2}w_1^2w_2\bigg(-\frac{\partial w_1}{\partial z_1}\bigg)+O\Big(|\beta|^2l^{2+\alpha}+|\beta|e^{-2l\sin\frac{\pi}{k}}\Big)\notag\\
=&-6k\int_{\mathbb{R}^2}w^2 (|z-\xi_1|)w^\prime(|z-\xi_1|)w(|z-\xi_2|)\frac{z_1-l}{|z-\xi_1|}\,dz+O\Big(|\beta|^2l^{2+\alpha}+|\beta|e^{-2l\sin\frac{\pi}{k}}\Big)\notag\\
=&-6k\int_{\mathbb{R}^2}w^2(|z|)w^\prime(|z|)\,w(|z+\xi_1-\xi_2|)\frac{z_1}{|z|}
\,dz\,+O\Big(|\beta|^2l^{2+\alpha}+|\beta|e^{-2l\sin\frac{\pi}{k}}\Big)\notag\\
=&c_0\,\cdot w(|\xi_1-\xi_2|)+O\Big(|\beta|^2l^{2+\alpha}+|\beta|e^{-2l\sin\frac{\pi}{k}}\Big)\notag\\
=&c_0\,\bigg(2\sin\frac{\pi}{k}\bigg)^{-1/2}\cdot
l^{-1/2}e^{-2l\sin\frac{\pi}{k}}+O\Big(|\beta|^2l^{2+\alpha}+|\beta|e^{-2l\sin\frac{\pi}{k}}\Big),
\end{align}where $c_0$ is a positive constant. Here we have used
the fact that (\ref{id1_7}) holds.

Combining (\ref{id5_4}) and (\ref{id5_5}), we see that
\begin{equation}\notag
c(l)\approx c_1l^{-1/2}e^{-2l\sin\frac{\pi}{k}}-c_2\beta l^{-3}+O\big(\beta l^{-5}\big).
\end{equation}
Moreover, we may choose $\hat{l}_\beta$ such that
\begin{equation}\label{id5_6}
\hat{l}_\beta^{-1/2}e^{-2\hat{l}_\beta\sin\frac{\pi}{k}}=\beta\hat{l}_\beta^{-3}.
\end{equation}
Then
\begin{equation}\notag
\hat{l}_\beta=\frac{1}{2\sin\frac{\pi}{k}}\log\frac{1}{\beta}+c_3\log\log\frac{1}{\beta}+c_4.
\end{equation}
Now we want to choose $l$ such that
$l\in(\hat{l}_\beta-\gamma,\hat{l}_\beta+\gamma)$ and $c(l)=0$ for
some suitable $\gamma>0$. It is remarkable that
$$
c\big(\hat{l}_\beta-\gamma\big)\approx
c_1\big(\hat{l}_\beta-\gamma\big)^{-1/2}e^{-2\hat{l}_\beta\sin\frac{\pi}{k}
+2\gamma\sin\frac{\pi}{k}} -c_2\beta\big(\hat{l}_\beta-\gamma\big)^{-3}+O\bigg(\beta\Big(\log\frac{1}{\beta}\Big)^{-5}\bigg)\notag\\
>0\,, $$ if $\gamma$ is large enough. Similarly,  $c\big(\hat{l}_\beta+\gamma\big)<0$
if $\gamma$ is large enough. Since $c(l)$ is continuous in $l$,
then by the mean-value theorem, there exists
$l_\beta\in(\hat{l}_\beta-\gamma,\hat{l}_\beta+\gamma)$ such that
$c(l_\beta)=0$. Consequently, the function $ \left(
\begin{array}{*{20}c}
    u_{l_\beta}+\phi_{l_\beta} \\ v_{d} e^{\sqrt{-1} \psi_{l_\beta}}
\end{array} \right)
=:
\left( \begin{array}{*{20}c}
    u_\beta \\ v_\beta
\end{array} \right)
$ is a solution of (\ref{id1_5}). Furthermore, it is easy to check
that $(u_\beta,v_\beta)$ satisfies all the properties of
Theorem~\ref{thm1_1}. Therefore we may complete the proof of
Theorem~\ref{thm1_1}.

\section{Proof of Theorem \ref{thm1_2}}

We first consider problem (\ref{ra2}) on a ball $B_R$:
\begin{equation}
\label{eqn-Br}
\left\{\begin{array}{l}
\Delta u - u+ u^3 +\beta u S^2=0, u=u(r), r<R \\
\Delta S -\frac{d^2}{r^2} S +S(1-S^2) +\beta u^2 S=0, S=S(r), r<R \\
 S(0)=0, S(R)=1, u(R)=0, u>0, 0<S<1
\end{array}
\right.
\end{equation}
Our idea is to find a solution of (\ref{eqn-Br}), and then let $R
\to +\infty$. To this end, we consider the associated energy
functional
\begin{eqnarray}
\label{ER} E_R[u, S] &=& \frac{1}{2} \int_{B_R} (|\nabla u|^2 +
u^2) +\frac{1}{2} \int_{B_R} (|\nabla S|^2 +\frac{d^2}{r^2} S^2)
\\ & & +\frac{1}{4} \int_{B_R} (1-S^2)^2 -\frac{\beta}{2}
\int_{B_R} S^2 u^2 -\frac{1}{4} \int_{B_R} u^4\,, \nonumber
\end{eqnarray}
for $u \in H_0^1 (B_R)$ and $S \in I_R\equiv\left\{ S\in H^1(B_R):
S(z)=S(|z|), S(0)=0, S(R)=1\right\}$.

Let the Nehari manifold be
\begin{equation*}
{\mathcal N}= \left\{ (u, S) \in H_0^1 (B_R) \times I_R,u \geq 0, u \not \equiv 0:
\int_{B_R} (|\nabla u|^2 +u^2) = \int_{B_R} (u^4 +\beta u^2 S^2)
\right\}\,.
\end{equation*}

Then we consider the following energy minimization problem
\begin{equation}\label{6.4}
c_R=\inf_{ (u, S)\in {\mathcal N} } E_R [u, S]\,,
\end{equation}
and we have
\begin{lem}\label{6.1}
If $\beta <0$, then $c_R$ is obtained by some radially symmetric
function $ (u_R, S_R)$. Furthermore, $u_R^{'}(r) <0$ and $S_R^{'}
(r)>0$ for $r>0$.
\end{lem}

\noindent {\bf Proof:} We follow the proof of (1) Theorem 3.3 of
\cite{LW1}. To this end, we define another energy functional
\begin{eqnarray}
\label{ER1} E_R^{\prime}[u, S] &=& \frac{1}{4} \int_{B_R} (|\nabla u|^2 +
u^2) +\frac{1}{2} \int_{B_R} (|\nabla S|^2 +\frac{d^2}{r^2} S^2)
\\ & & +\frac{1}{4} \int_{B_R} (1-S^2)^2 -\frac{\beta}{4}
\int_{B_R} S^2 u^2\, \nonumber
\end{eqnarray}
and another solution manifold
\begin{equation*}
{\mathcal N}^{\prime}= \left\{ (u, S) \in H_0^1 (B_R) \times I_R, u \geq 0, u \not \equiv 0 :
\int_{B_R} (|\nabla u|^2 +u^2)  \leq \int_{B_R} (u^4 +\beta u^2 S^2)
\right\}\,.
\end{equation*}

We consider another minimization problem:
\begin{equation}\label{6.4n}
c_R^{\prime}=\inf_{ (u, S)\in {\mathcal N}^{\prime} } E_R^{\prime} [u, S].
\end{equation}

Certainly, we have
\begin{equation}
\label{6.5n}
c_R^{'} \leq c_R.
\end{equation}

Let $ (u_n, S_n)$ be a minimizing sequence of $c_R^{\prime}$ on
${\mathcal N}^{\prime}$. Replacing $S_n$ by $ \min(S_n, 1)$, we may assume
that $S_n \leq 1$. We may denote $u_n^*$ and $S_n^*$ as the
Schwartz symmetrization of $u_n$ and $S_n$, respectively. Then
$(1-S_n)^{*}= 1-S_n^{*}$. By Theorem~3.4 of~\cite{lb},
\begin{equation}
\int_{B_R} (1-S_n^2) u_n^2 \leq \int_{B_R} (1-S_n^2)^{*} (u_n^{*})^2
\end{equation}
and hence due to $\beta <0$,
\begin{equation}
-\beta \int_{B_R} (u_n^{*})^2 (S_n^{*})^2 \leq -\beta \int_{B_R} u_n^2 S_n^2
\end{equation}

On the other hand, we also have
\[ \int_{B_R} (\frac{1}{2} (u_n^{*})^2 + \frac{1}{2} \frac{d^2}{r^2} ( S_n^{*})^2 -\frac{1}{4} (u_n^{*})^4)
= \int_{B_R} (\frac{1}{2}  u_n^2 + \frac{1}{2} \frac{d^2}{r^2}  S_n^2 -\frac{1}{4} u_n^4)\,, \]
\[ \int_{B_R} (|\nabla u_n^{*}|^2 +|\nabla S_n^{*}|^2) \leq \int_{B_R} (|\nabla u_n|^2 +|\nabla S_n|^2).\]
Hence we obtain
\begin{equation*}
E_R^{'} [u_n^{*}, S_n^{*}]\leq E_R^{'} [u_n, S_n]
\end{equation*}
and
\begin{equation}
\label{C1}
\int_{B_R} (|\nabla u_n^{*}|^2 +(u_n^{*})^2) \leq  \int_{B_R} ((u_n^{*})^4 +\beta (u_n^{*})^2 (S_n^{*})^2).
\end{equation}
Thus we may replace $ (u_n, S_n)$ by its symmetrization $(u_n^{*},
S_n^{*})$. Since $ S_n \leq 1$ and $ H^1 (B_R)$ is a compact
embedding to $L^4 (B_R)$, we see that $(u_n S_n) \to (u_R, S_R) $
 weakly in $H^1 (B_R)$ and strongly in $L^4 (B_R)$, where $(u_R, S_R) \in {\mathcal N}^{\prime}$ attains $c_R^{'}$. So $c_R^{'}$ is attained. If $(u_R, S_R) \in ({\mathcal N}^{\prime})^{\circ}-$the
interior of ${\mathcal N}^{\prime}$, then $(u_R, S_R)$ is a local minimizer
of $E_R^{\prime}$ and hence we have
\begin{equation*}
\Delta u_R -u_R +\beta u_R S_R=0, \quad \forall 0<r<R, \quad\hbox{
and }\: u_R (R)=0\,,
\end{equation*}
which implies $ u_R \equiv 0$ since $\beta <0$. This is impossible
since from (\ref{C1}) and Sobolev embedding, we infer that
 $ \int_{B_R} u_R^4 \geq C>0$. Therefore $(u_R, S_R) \in \partial ({\mathcal N}^{\prime})= {\mathcal N}$ and hence
\begin{equation}
\label{6.6n}
 c_R \leq E_R [u_R, S_R]= E_R^{'} [u_R, S_R] = c_R^{'}.
\end{equation}

Combining (\ref{6.5n}) and (\ref{6.6n}), we conclude that $c_R$ is attained by $(u_R, S_R)$.  Then we have the following equality:
\begin{equation}\label{6.1.0} G_R [u_R, S_R]=\int_{B_R} (|\nabla u_R|^2 + u_R^2
-\beta u_R^2 S_R^2 -u_R^4)=0\,.\end{equation} Hence there exists a
Lagrange multiplier $\lambda_R$ such that
\begin{equation}\label{6.1.1} \nabla E_R +\lambda_R \nabla G_R=0\,.\end{equation}
Acting (\ref{6.1.1}) on $(u_R,0)$, we may obtain
\[\int_{B_R} (|\nabla u_R|^2 +u_R^2 -\beta u_R S_R^2 -u_R^4)
+2\lambda_R \int_{B_R} (|\nabla u_R|^2 +u_R^2 -\beta u_R^2 S_R^2
-2 u_R^4)=0\,,\] and hence by (\ref{6.1.0}),
\[ -\lambda_R\,\int_{B_R} u_R^4 =0\,,\quad\hbox{ i.e. }\quad\lambda_R=0\,.\]
Therefore, we may complete the proof of Lemma~\ref{6.1}.

\qed

Theorem~\ref{thm1_2} is proved by the following lemma

\begin{lem}
As $R \to +\infty$, $ (u_R, S_R) \to (u_\infty, S_\infty)$ and $(u_\infty, S_\infty)$ is a solution of (\ref{ra2}).
\end{lem}

\noindent {\bf Proof:} Since $S_R \leq 1, S_R (0)=0$, we first
show that $u_R$ is uniformly bounded, independent of $R>1$.
Actually, it is sufficient to show that $\int_{B_R} (|\nabla
u_R|^2+u_R^2) \leq C$, where $C$ is a positive constant
independent of $R>1$. Let
\[ GL_{B_R} [S]= \int_{B_R} \left( \frac{1}{2} |\nabla S|^2 +\frac{d^2}{2r^2} S^2 +\frac{1}{4}(1-S^2)^2\right)\,.\]
It is remarkable that $GL_{B_{R}}$ may come from the conventional
Ginzburg-Landau functional.

From~\cite{PR}, we may set $\bar{S}_R$ as the unique minimizer of
$GL_{R}$. Then for any $ u \in H_0^1 (B_R)$, there exists $t_R$
such that $ (\sqrt{t_R} u, \bar{S}_R) \in {\mathcal{N}}$,
where $t_R$ is simply given by
\begin{equation}
 t_R = \frac{\int_{B_R} (|\nabla u|^2 +u^2 -\beta \bar{S}_R^2 u^2)}{\int_{B_R} u^4
 }\,.
\end{equation}
Thus by  Lemma~\ref{6.1} and $0\leq \bar{S}_R<1$,
\begin{eqnarray*}
 c_{R} \leq E_R[\sqrt{t_R} u, \bar{S}_R]
 &=& GL_{B_R} (\bar{S}_R) + \frac{1}{4}\,\left[\frac{\int_{B_R} |\nabla u|^2 + u^2 -\beta \bar{S}_R^2 u^2}{ (\int_{B_R}
 u^4)^{\frac{1}{2}}}\right]^2 \\ &\leq& GL_{B_R} (\bar{S}_R) + \left[\frac{\int_{B_R} |\nabla u|^2 + (1-\beta) u^2}{ (\int_{B_R}
 u^4)^{\frac{1}{2}}}\right]^2\,,
\end{eqnarray*} for all $ u \in H_0^1 (B_R)$. Consequently, due to
$\beta<0$,
\begin{equation}\label{upper} c_R\leq GL_{B_R} (\bar{S}_R) +
C_0\,,\end{equation} where $C_0$ is a positive constant
independent of $R>1$. Here we have used the fact that
$\ds\lim_{R\to\infty}\inf_{u\in H^1_0(B_R)}\,\frac{\int_{B_R}
|\nabla u|^2 + (1-\beta) u^2}{ (\int_{B_R}
u^4)^{\frac{1}{2}}}<\infty\,.$

By standard theory of Ginzburg-Landau equation (cf.~\cite{BBH}),
we have
\begin{equation}\label{lower}
GL_{B_R} [S_R]\geq GL_{B_R} [\bar{S}_R]\,.
\end{equation}
Combining (\ref{upper}) and (\ref{lower}), we see that
\[  \frac{1}{2} \int_{B_R} (|\nabla u_R|^2 + u_R^2 -\beta S_R^2 u_R^2) -\frac{1}{4}\int_{B_R} u_R^4 \leq C_3\]
and hence by the equation of $u_R$,
\begin{equation}
\int_{B_R} (|\nabla u_R|^2 + u_R^2) \leq C\,,
\end{equation}
from which standard elliptic regularity theory gives that $u_R
\leq C$. Thus we may obtain that $(u_R, S_R) \to (u_\infty,
S_\infty)$ which solves $\Delta u_\infty -u_\infty + u_\infty^3
+\beta u_\infty S_\infty^2=0$. Note that $ u_R(0)\geq 1$ and hence
$u_\infty \not \equiv 0$. By the Maximum Principle, $u_\infty >0$.
Similarly, $ 0<S_\infty(r) <1 $ for $r>0$. Therefore we may
complete the proof of Theorem~\ref{thm1_2}.

\qed

\section{Proof of Theorem \ref{thm1_3}}
\noindent

In this section, we want to construct $S^2$-valued map to get
half-skyrmions by spike-vortex solutions obtained in
Theorem~\ref{thm1_1} and \ref{thm1_2}. For simplicity, we firstly
use spike-vortex solutions in Theorem~\ref{thm1_2} to find
half-skyrmions. Let $ (u, v)$ be the radial spike-vortex solution
in Theorem~\ref{thm1_2}. We may define a $S^2$-valued map by
\beq\label{sk1} \overrightarrow{n}\equiv\frac{1}{\sqrt{u^2+v_1^2+v_2^2}}\,\left(\begin{array}{ll}&v_1\\
&v_2\\ &u\end{array}\right)=\left(\begin{array}{ll}
&\cos(\phi(r))\,\cos(d\theta) \\ &\cos(\phi(r))\,\sin(d\theta) \\
&\sin(\phi(r))\end{array}\right)\,,\eeq where
$v=v_1+\sqrt{-1}\,v_2$, \beq\label{sk2} \cos(\phi(r))=
\frac{f(r)}{\sqrt{u^2+f^2}}\,, \eeq and \beq\label{sk3}
\sin(\phi(r))=\frac{u(r)}{\sqrt{u^2+f^2}} \,.\eeq Since both $u$
and $f$ are positive everywhere, the function $\phi$ is
well-defined and single-valued. The map $\overrightarrow{n}$ can
be decomposed into
$$ \overrightarrow{n}= \cos(\phi(r))\left(\begin{array}{ll}
&\cos(d\theta) \\ &\sin(d\theta) \\
&0\end{array}\right) + \sin(\phi(r))\left(\begin{array}{ll}&0\\
&0\\ &1\end{array}\right)\,.$$ Then it is easy to check that
\begin{eqnarray*} \int_{\mathbb{R}^2}\overrightarrow{n}\cdot\left(\pp_x
\overrightarrow{n}\wedge\pp_y\overrightarrow{n}\right)\,dx\,dy &=&
\int_{\mathbb{R}^2}\,-\frac{d}{r}\phi^\prime(r)\cos(\phi(r))\,dx\,dy
\\ &=& -2\pi\,d \sin(\phi(r))|_{r=0}^{+\infty} =
2\pi\,d\,,
\end{eqnarray*} i.e.
 the topological charge \beq\label{sk4}
Q=\frac{1}{4\pi}\,\int_{\mathbb{R}^2}\overrightarrow{n}\cdot\left(\pp_x
\overrightarrow{n}\wedge\pp_y\overrightarrow{n}\right)\,dx\,dy=
\frac{d}{2}\,.\eeq Here we have used (\ref{sk3}) and the fact that
$u(0)>0$, $u(+\infty)=f(0)=0$ and $f(+\infty)=1$.

For the spike-vortex solution $(u,v)$ in Theorem~\ref{thm1_1},
since $\beta>0$ sufficiently small, the associated map
$\overrightarrow{n}$ has the following form \beq\label{sk11}
\overrightarrow{n}\equiv\frac{1}{\sqrt{u^2+v_1^2+v_2^2}}\,\left(\begin{array}{ll}&v_1\\
&v_2\\ &u\end{array}\right)=\left(\begin{array}{ll}
&\cos(\phi)\,\cos(d\psi) \\ &\cos(\phi )\,\sin(d\psi) \\
&\sin(\phi)\end{array}\right)\,,\eeq  where $\phi=\phi(r,\theta)$
and $\psi=\psi(r,\theta)$ satisfying \beq\label{sk2.1}
\cos(\phi(r,\theta))= \frac{|v|(r,\theta)}{\sqrt{u^2+|v|^2}}\,,
\eeq \beq\label{sk3.1}
\sin(\phi(r,\theta))=\frac{u(r,\theta)}{\sqrt{u^2+|v|^2}} \,,\eeq
and $\psi=\theta+h$, where $h$ is a single-valued regular function
satisfying $h=O(\beta)$ as $\beta\to 0+$. Here both $u$ and $|v|$
may not have radial symmetry. Due to $\beta>0$, we may apply the
standard maximum principle on the first equation of the system
(\ref{id1_5}) i.e. the equation of $u$. Then the solution $u$ is
positive everywhere so the function $\phi$ is well-defined and
single-valued.

Now we want to calculate the topological charge $Q$ as for
(\ref{sk4}). By~(\ref{sk11}), it is easy to check that
\begin{equation*} \overrightarrow{n}\cdot\left(\pp_x
\overrightarrow{n}\wedge\pp_y\overrightarrow{n}\right)=
-\frac{d}{r}\,\phi_r\cos\phi + \frac{d}{r}\left(\phi_\theta h_r
-\phi_r h_\theta\right)\cos\phi\,.
\end{equation*} Hence by (\ref{sk2.1}), (\ref{sk3.1}) and using integration by part, we may obtain
\begin{equation*} Q=\frac{1}{4\pi}\,\int_{\mathbb{R}^2}\overrightarrow{n}\cdot\left(\pp_x
\overrightarrow{n}\wedge\pp_y\overrightarrow{n}\right)=
\frac{d}{2}\,.
\end{equation*} Here we have used the fact that $u(0)>0$, $u(\infty)=0$,
$v(0)=0$ and $|v(\infty)|=1$. Therefore we may complete the proof
of Theorem~\ref{thm1_3}.

\end{document}